\documentclass{IEEEtran4PSCC}
\ifCLASSINFOpdf
   \usepackage[pdftex]{graphicx}
\else
   \usepackage[dvips]{graphicx}
\fi
\usepackage[cmex10]{amsmath}
\usepackage{algorithm}
\usepackage{algpseudocode}
\usepackage{amssymb}
\usepackage{tabularx}
\usepackage{siunitx} 
\usepackage{caption}
\usepackage{subcaption}
\usepackage{multirow}
\usepackage{booktabs}
\usepackage{tikz}
\usepackage{amsfonts}
\usepackage{amsmath}
\usepackage{svg}
\usepackage{fontawesome5} 
\usepackage{graphicx,svg,circuitikz}
\usetikzlibrary{positioning, shapes, arrows.meta, calc, fit, backgrounds}
\usetikzlibrary{shapes.geometric, arrows,positioning,arrows.meta,calc,fit,backgrounds, shapes}
\definecolor{tugreen}{RGB}{132, 184, 25}
\definecolor{lightgray}{RGB}{220,220,220}
\definecolor{titlegreen}{RGB}{132,184,25}
\definecolor{customorange}{HTML}{CA7406}
\usetikzlibrary{decorations.pathreplacing}
\usepackage{xcolor}
\usepackage{soul}
\setstcolor{red}

\usepackage{tikz-3dplot}
\usepackage{etoolbox} 
\robustify\bfseries
\makeatletter
\let\old@ps@headings\ps@headings
\let\old@ps@IEEEtitlepagestyle\ps@IEEEtitlepagestyle
\def\psccfooter#1{%
    \def\ps@headings{%
        \old@ps@headings%
        \def\@oddfoot{\strut\hfill#1\hfill\strut}%
        \def\@evenfoot{\strut\hfill#1\hfill\strut}%
    }%
    \def\ps@IEEEtitlepagestyle{%
        \old@ps@IEEEtitlepagestyle%
        \def\@oddfoot{\strut\hfill#1\hfill\strut}%
        \def\@evenfoot{\strut\hfill#1\hfill\strut}%
    }%
    \ps@headings%
}
\makeatother

\psccfooter{%
        \parbox{\textwidth}{\hrulefill \\ \small{24th Power Systems Computation Conference} \hfill \begin{minipage}{0.2\textwidth}\centering \vspace*{4pt} \includegraphics[scale=0.06]{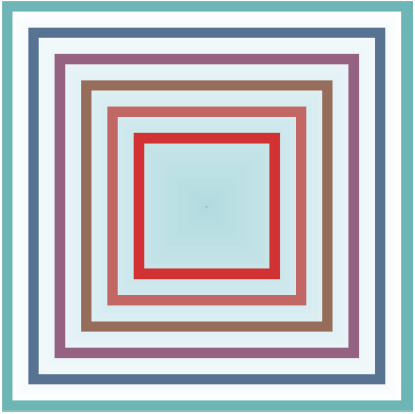}\\\small{PSCC 2026} \end{minipage} \hfill \small{Limassol, Cyprus --- June 8-12, 2026}}%
}

\usepackage[style=ieee, url=false,doi=false,isbn=false,date=year,citestyle=numeric-comp,maxbibnames=10,minbibnames=10,maxcitenames=10,mincitenames=10]{biblatex}
\bibliography{lit_clean.bib,MBB_lit.bib}

\begin{document}
%
\title{Towards grid-aware multi-period flexibility aggregation --  A constrained zonotope approach}


\author{\IEEEauthorblockN{Maurice Raetsch\IEEEauthorrefmark{1},
Maísa Beraldo Bandeira\IEEEauthorrefmark{1},
Christian Rehtanz\IEEEauthorrefmark{1}, 
Alexander Engelmann\IEEEauthorrefmark{2} and
Timm Faulwasser\IEEEauthorrefmark{3}}
\IEEEauthorblockA{\IEEEauthorrefmark{1} ie3\\
TU Dortmund University,
Dortmund, Germany\\ \{maurice.raetsch, maisa.bandeira, christian.rehtanz\}@tu-dortmund.de}
\IEEEauthorblockA{\IEEEauthorrefmark{2}
logarithmo GmbH \& Co. KG\\
Dortmund, Germany\\ alexander.engelmann@ieee.org}
\IEEEauthorblockA{\IEEEauthorrefmark{3} Institute of Control Systems\\
Hamburg University of Technology,
Hamburg, Germany\\ timm.faulwasser@ieee.org}
}

\maketitle

\begin{abstract}
Aggregation schemes provide a means to reduce the computational complexity of power system operation by reducing the number of devices that are considered individually. This can be achieved with tools of computational geometry, where the feasible set is projected onto the decision variables of the point of interconnection. Set projection is computationally expensive, especially in the context of multi-period power system operation. This calls for efficiency improvements via structure exploitation of set representations. This paper proposes efficient flexibility aggregation via constrained zonotopes. We evaluate the performance of the proposed method on a 15-bus distribution grid with time-dependent elements for up to 96 timesteps. The results suggest that the presented method significantly improves computation times compared to classic polytope projection approaches.
\end{abstract}


\thanksto{\noindent TF and CR acknowledge financial support from TRR 391 Spatio-temporal Statistics for the
Transition of Energy and Transport (520388526) by the Deutsche Forschungsgemeinschaft (DFG,
German Research Foundation).\\\noindent Submitted to the 24th Power Systems Computation Conference (PSCC 2026).}

\section{Introduction}\label{sec:intro}
The increasing integration of Distributed Energy Resources (DERs) into modern power systems introduces new operational challenges for both Transmission System Operators (TSOs) and Distribution System Operators (DSOs). From the perspective of the TSO, higher volatility and the spatial distribution of generation increase the risk of voltage violations, line congestion, and inverted power flows. Moreover, most flexible devices are located in distribution grids added at the distribution level. The availability of DERs and storage solutions in the distribution grid offers new opportunities for flexible grid operation, which in turn creates new coordination requirements between the TSO and the DSOs to ensure reliable and efficient system performance.

These interactions can be modeled through a centralized Optimal Power Flow (OPF) formulation, in which all system data—such as DER capacities, grid parameters, and asset locations—are collected by a central operator that solves a single global optimization problem. However, to preserve data privacy, reduce computational complexity and to simplify information exchange, such detailed information from the distribution grids should not be directly shared with the TSO.



To address these challenges, hierarchical optimization has been widely adopted in the literature. In the first stage of the approach, flexibility aggregation, the Feasible Operation Region (FOR) is computed. The FOR describes the values of the shared variables, usually active and reactive power, at the interconnection point that can be achieved without compromising the feasibility of the optimization problem from the DSO's perspective. A survey of different methods can be found in \cite{sarstedt2021}. The TSO can then use this information to determine the optimal utilization of DER flexibility without requiring direct data from the distributed grid assets.

The standard FOR estimation methods rely on sampling the region and computing the convex hull of these points. The samples can be randomly obtained \cite{Gonzales2018} or calculated in a more structured manner, utilizing knowledge of the grid and constraints to solve optimization problems \cite{silva2018, sarstedt2022, capitanescu2018}. The methods above consider assumptions that can be quite limiting. They assume a unique interconnection point between DSO and TSO, the voltage at this point to be constant, and no time-dependent elements, e.g., storage systems. To mitigate some of these limitations, recent optimization-bases approaches compute inner-approximations of the aggregated flexibility in multi-period settings \cite{Chen2020,Chen2021}. Similar approaches extend this approach to consider uncertainties \cite{Cui2021} or solve optimization problems to parameterize a virtual battery model considering time-dependent elements and network constraints \cite{2Tan2024}, but these inner approximation based methods remain inherently constrained by the geometry of the chosen surrogates.


There is limited research on the flexibility aggregation in distribution grids that considers time-dependent elements. A survey of various methods for aggregating batteries can be found in \cite{ozturk2022,ozturk2024alleviating}. They offer computational tractability but neglect grid constraints or rely on data-driven formulations, which reduces their accuracy for grid-aware applications. Multi-period cost curve formulations \cite{MP_CC_Aggregation} capture temporal coupling but limits the FOR to only portraying active power. 




Recently, set projection has been used for the computation of the FOR~\cite{2Tan2024, Engelmann2025}. The method has been extended to FORs, which include time-dependent elements~\cite{wang2023}, but it requires solving multiple optimization problems. 
For affine grid models, such as LinDistFlow \cite{baran1989-0}, the projection can be obtained in one step via tools of computational geometry. The downside of this technique is the computational burden of the set projection operation, which limits its scalability \cite{Bandeira2025}.

One line of research addresses this computational burden of polytope projection by applying zonotope approximations~\cite{agg_disagg_zonotope,zono_union}. However, such approximations are not always accurate, particularly when the underlying set is non-symmetric.

The present paper introduces a novel approach to compute time-coupled and voltage dependent active and reactive power FORs by extending the non-iterative projection-based formulation~\cite{ADP} to multi-period settings with time-dependent elements. To overcome computational bottlenecks that normally hinder the scalability of set projections, we rely on Constrained Zonotopes ($\mathcal{CZ}$)~\cite{CZ}. $\mathcal{CZ}$s enable non-symmetric set representation with efficient set projection properties. This results in significant reduction of computation times without compromising representation accuracy. Furthermore, most of the remaining computation can be carried out offline, enabling the practical implementation for the computation of FORs for a 15-bus distribution grid that include time-dependent elements such as storage units for up to 96 timesteps.




The remainder of this paper is structured as follows: Section~\ref{sec:probl_statement} elaborates the problem statement and recaps the flexibility aggregation via polytope projection. In Section~\ref{sec:meth_cZ} we present how $\mathcal{CZ}\mathrm{s}$ can be utilized for efficient flexibility aggregation over time, while Section~\ref{sec:additional_intersec} highlights the potential of this method for dynamic problem settings. Section~\ref{sec:numerical_case_study} compares and discusses numerical results for the presented methods. The paper concludes with Section~\ref{sec:conclusion}.

\section{Problem Statement} \label{sec:probl_statement}
\begin{figure}[t!]
    \centering
    \includegraphics[scale=0.6]{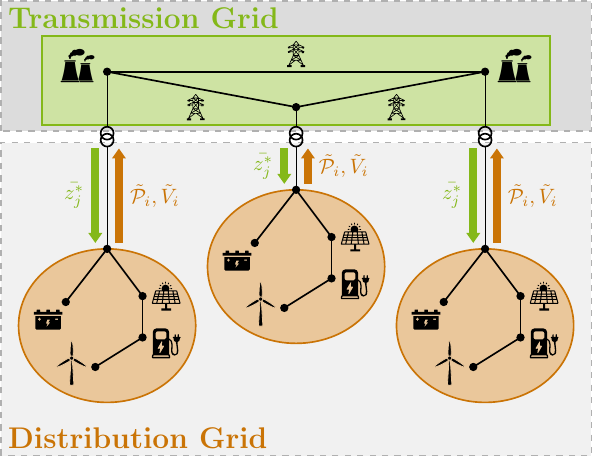}
    \caption{Tree-structured optimization of power system operation adapted from \cite{ADP}.}
    \label{fig:example_adp}
\end{figure}

We assume a tree-structured power system, divided into transmission and distribution grids, as shown in Figure \ref{fig:example_adp}. Time coupled system operation across multiple grid levels requires the solution of multi-stage Optimal Power Flow (OPF) problems. Next, we introduce one variant of this problem, which we later use as a basis for aggregation. We limit the presentation to one distribution grid in order to simplify presentation.


The first step in flexibility aggregation via set projection is to define the constraint set of the sub-problems, i.e. the distribution grids. We represent the distribution grid as a graph  $G^e = (\mathcal N_d, \mathcal B_d)$, where $\mathcal N_d$ is the set of all buses and $\mathcal B_d \subseteq \mathcal N_d \times \mathcal N_d$ is the set of branches, i.e., lines and transformers. 

We define the nodal active and reactive power $p_{m}(k)$ and $q_{m}(k)$ of node $m$ at each timestep $k$ in the considered interval $\mathcal I = {1,\dots,N}$ as \begin{align} \label{eq:PFeq}
	&p_m(k) = \sum_{l \in \mathcal{N}_{d}} p_{m,l}(k) , \qquad  q_m = \sum_{l \in \mathcal{N}_{d}} q_{m,l}(k),
\end{align} 
with $p_{m,l}(k)$ and $q_{m,l}(k)$ representing the active and reactive power flow between nodes $m$ and $l$ at each timestep $k$. The net powers at each node are 
\begin{subequations}\label{eq:netPwr}
    \begin{align} 
	p_m(k) &= p_m^g(k)   - p_m^d(k) + p_m^s(k), \\ 	q_m(k) &= q_m^g(k) - q_m^d(k),
    \end{align}
\end{subequations} where $p_m^g(k)$ and $q_m^g(k)$ are the active/reactive power generation and $p_m^d$ and $q_m^d(k)$ are the active/reactive power demand at node $m \in \mathcal N_d$. Furthermore, $p_m^s(k)$ is the battery's active power. For nodes without a battery, $p_m^s = 0$, while for nodes with a battery, $p_m^s$ is negative when charging and positive when discharging, i.e., providing the grid with power.

The active power generation for the flexible renewable generators must stay within the 
bounds
\begin{align}
	 0 \leq p_m^g(k) \leq \bar f_m^g(k), \qquad m \in \mathcal G. \label{eq:genCstrBds}
\end{align} 
where $\mathcal{G}\subseteq \mathcal{N}$ denotes the set of buses with controllable generators and $f_m^g(k)$ represents the maximal active power generation given by irradiation/wind forecast conditions at each timestep.

We limit the reactive power support from each renewable generator by constraining the maximum apparent power $\bar s_m(k)$ and the power factor limit $\cos(\phi)$ at each timestep with the affine constraints \cite{contreras2019}
\begin{subequations} \label{eq:ice_cream_ldf}
\begin{align}
    p_m^g(k) &\leq \bar s_m(k)\cos(\phi),\\
    -p_m^{g}(k)\tan(\phi) &\leq q_m^{g}(k)\leq p_{m}^{g}(k)\tan(\phi).
\end{align}
\end{subequations} 

Next, we formulate the active/reactive power flow $p_{m,l}(k)$ and $q_{m,l}(k)$ on each line $(m,l) \in \mathcal B_d$ and the squared voltage $\nu_m$ at each bus $m \in \mathcal N_d$  according to the DistFlow formulation for radial grids \cite{baran1989}
\begin{subequations} \label{eq:PF_df} \begin{align} 
	p_{m,l}(k)  &= \sum_{j \in \mathcal N_d} p_{m,j}(k) + p_m(k) + r_{m,l}\,\ell_{m,l}(k),\label{eq:df_p} \\ 
	  q_{m,l}(k) & = \sum_{j \in \mathcal N_d}   q_{m,j}(k) + q_m(k)+x_{m,l}\ell_{m,l}(k) ,\label{eq:df_q} \\
	  \nu_m(k) &= \nu_l(k) + 2(r_{m,l}\,p_{m,l}(k) + x_{m,l}\,q_{m,l}(k)) \notag \\&-(r_{m,l}^2+x_{m,l}^2) \ell_{m,l}(k),\label{eq:df_v}
\end{align}  \end{subequations} 
where, $r_{m,l}$ and $x_{m,l}$ are, respectively, the line resistance and reactance and $\ell_{m,l}(k)$ is the squared line current.

To obtain a polyhedral constraint set, it is common to neglect line impedances. Then the squared
branch current $\ell_{m,l}(k)$ does not need to be calculated. As shown in \cite{ADP}, this results in a relatively good approximation of the true flexibility region with a slight offset due to the neglected losses. Similar to \cite{nazir2022,jha2020,Brett2016} we use the first-order Taylor series as linear approximation for the square current $\ell_{m,l}(k)$ around the current operational point $t_{m,l}^0$ \cite{jha2020},
\begin{align}\label{eq:approximated_l}
 	\ell_{m,l}(k) \approx \ell_{m,l}(k)^0 + J_{m,l}^{\top}|_{\substack{t_{m,l}^0(k)}} \delta_{m,l}(k) 
\end{align}
where 
\begin{align} \label{eq:delta}
 	\delta_{m,l}(k) &= \begin{bmatrix} p_{m,l}(k)-p_{m,l}^0 \\ q_{m,l}(k)-q_{m,l}^0 \\ \nu_{l}(k)-\nu_{l}^0 \end{bmatrix},\\
    J_{m,l}|_{\substack{t_{m,l}^0}} &=  \begin{bmatrix} \frac{2p_{m,l}^{0}}{\nu_l^0} \\ 
 		\frac{2q_{m,l}^0}{\nu_l^{0}}  \\
 		- \frac{(p_{m,l}^0)^2+(q_{m,l}^0)^2 }{(\nu_l^0)^2} \end{bmatrix}. 
\end{align}
Additionally, we limit squared voltages and currents 
\begin{align} 
		\underline \nu_m \leq \nu_m(k) \leq \bar \nu_m, \quad \underline \ell_{m,l} \leq \ell_{m,l}(k) \leq \bar \ell_{m,l},\label{eq:box_constraints}
\end{align} 

Let $\mathcal{S} \subseteq \mathcal{N}$ denote the set of buses equipped with storage units. Each storage unit has a limited capacity $e_m^s(k)$
\begin{equation}\label{eq:prel_storage_limits}
\underline{e}^s_m \leq e^s_m(k) \leq \overline{e}^s_m,\quad m\in\mathcal{S}
\end{equation}
with bounded charge/discharge rate
\begin{equation}\label{eq:prel_storage_power_limits}
\underline{p}^s_m \leq p_m^s(k) \leq \overline{p}^s_m,\quad m\in\mathcal{S}.
\end{equation}
 For the sake of simplicity, we neglect the reactive power support of   storage systems as well as the charging and discharging efficiency of the storage units.

A simple discrete-time integrator dynamic models the State of Charge (SOC) of the storage units. Assuming a fixed timestep  $\Delta k$ and neglecting losses, the SOC dynamic reads
\begin{equation}\label{eq:prel_storage_dynamics}
e_m(k+1) = e_m(k) - \Delta k \cdot p_m^s(k),\quad  e_m(1) = e_m^0
\end{equation}
with $m\in\mathcal{S}$.
We define the resulting constraint set as 
\begin{align*}
	\mathcal X_d \hspace{-.1cm}\doteq\hspace{-.1cm} \big \{ x_i \in \mathbb{R}^{n_{xi}} \,|\,&\eqref{eq:PFeq}\text{-}\eqref{eq:prel_storage_dynamics} \text{ hold } \forall m \in \mathcal N_d,\, \\ &\forall k \in \mathcal I,\,\text{and } \forall (m,l) \in \mathcal B_m\big \}.
\end{align*}

The vector of coupling variables with the TSO is defined as  $z_i^\top=[[p_{i,j}(k),q_{i,j}(k)]_{i,j\in \mathcal{B}_d, k \in \mathcal I},[\nu_i(k)]_{i\in\mathcal{N}_d^c, k \in \mathcal I}]$.
These variables appear in both the TSO and DSO optimization problems, representing the physical coupling between the two grids. Here, $\mathcal{N}_d^c$ denotes the set of nodes that serve as coupling points between the transmission and distribution grids.



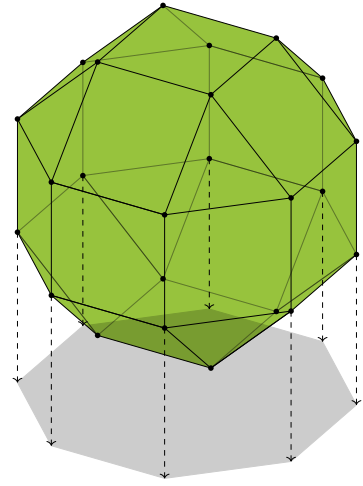
\begin{figure}[t!!]
    \centering
    \vspace{-30pt}

\tdplotsetmaincoords{60}{120}
\definecolor{green}{RGB}{132,184,25}
\definecolor{grey}{RGB}{0,0,0}
\begin{tikzpicture}[tdplot_main_coords, scale=3,vector guide/.style={dashed,black,->}]
\scalebox{0.7}{
\coordinate (v1)  at ( 0.414,  0.414,  1);
\coordinate (v2)  at ( 0.414, -0.414,  1);
\coordinate (v3)  at (-0.414,  0.414,  1);
\coordinate (v4)  at (-0.414, -0.414,  1);

\coordinate (v5)  at ( 0.414,  0.414, -1);
\coordinate (v6)  at ( 0.414, -0.414, -1);
\coordinate (v7)  at (-0.414,  0.414, -1);
\coordinate (v8)  at (-0.414, -0.414, -1);

\coordinate (v9)  at ( 0.414,  1,  0.414);
\coordinate (v10) at ( 0.414, -1,  0.414);
\coordinate (v11) at (-0.414,  1,  0.414);
\coordinate (v12) at (-0.414, -1,  0.414);

\coordinate (v13) at ( 0.414,  1, -0.414);
\coordinate (v14) at ( 0.414, -1, -0.414);
\coordinate (v15) at (-0.414,  1, -0.414);
\coordinate (v16) at (-0.414, -1, -0.414);

\coordinate (v17) at ( 1,  0.414,  0.414);
\coordinate (v18) at ( 1,  0.414, -0.414);
\coordinate (v19) at ( 1, -0.414,  0.414);
\coordinate (v20) at ( 1, -0.414, -0.414);

\coordinate (v21) at (-1,  0.414,  0.414);
\coordinate (v22) at (-1,  0.414, -0.414);
\coordinate (v23) at (-1, -0.414,  0.414);
\coordinate (v24) at (-1, -0.414, -0.414);

\draw[fill=green,opacity=0.6, draw=none] (v1)--(v2)--(v4)--(v3)--cycle;
\draw[fill=green,opacity=0.6, draw=none] (v5)--(v6)--(v8)--(v7)--cycle;

\draw[fill=green,opacity=0.6, draw=none] (v1)--(v9)--(v11)--(v3);
\draw[fill=green,opacity=0.6, draw=none] (v2)--(v10)--(v12)--(v4);

\draw[fill=green,opacity=0.6, draw=none] (v5)--(v13)--(v15)--(v7);
\draw[fill=green,opacity=0.6, draw=none] (v6)--(v14)--(v16)--(v8);

\draw[fill=green,opacity=0.6, draw=none] (v1)--(v17)--(v19)--(v2);
\draw[fill=green,opacity=0.6, draw=none] (v3)--(v21)--(v23)--(v4);

\draw[fill=green,opacity=0.6, draw=none] (v5)--(v18)--(v20)--(v6);
\draw[fill=green,opacity=0.6, draw=none] (v7)--(v22)--(v24)--(v8);

\draw[fill=green,opacity=0.6, draw=none] (v10)--(v14)--(v16)--(v12);
\draw[fill=green,opacity=0.6, draw=none] (v11)--(v15)--(v13)--(v9);

\draw[fill=green,opacity=0.6, draw=none] (v9)--(v17)--(v1);
\draw[fill=green,opacity=0.6, draw=none] (v10)--(v19)--(v2);
\draw[fill=green,opacity=0.6, draw=none] (v11)--(v21)--(v3);
\draw[fill=green,opacity=0.6, draw=none] (v12)--(v23)--(v4);

\draw[fill=green,opacity=0.6, draw=none] (v13)--(v18)--(v5);
\draw[fill=green,opacity=0.6, draw=none] (v14)--(v20)--(v6);
\draw[fill=green,opacity=0.6, draw=none] (v15)--(v22)--(v7);
\draw[fill=green,opacity=0.6, draw=none] (v16)--(v24)--(v8);

\draw[fill=green,opacity=0.6, draw=none] (v17)--(v18)--(v13)--(v9);
\draw[fill=green,opacity=0.6, draw=none] (v19)--(v20)--(v18)--(v17);
\draw[fill=green,opacity=0.6, draw=none] (v22)--(v24)--(v23)--(v21);
\draw[fill=green,opacity=0.6, draw=none] (v23)--(v24)--(v16)--(v12);
\draw[fill=green,opacity=0.6, draw=none] (v21)--(v22)--(v15)--(v11);
\draw[fill=green,opacity=0.6, draw=none] (v10)--(v14)--(v20)--(v19);

\draw (v1)--(v2); \draw (v2)--(v4); \draw (v4)--(v3); \draw (v3)--(v1);
\draw (v5)--(v6); \draw[opacity=0.4] (v6)--(v8); \draw[opacity=0.4] (v8)--(v7); \draw[opacity=0.4] (v7)--(v5);

\draw (v1)--(v9); \draw (v9)--(v11); \draw (v11)--(v3);
\draw (v2)--(v10); \draw (v10)--(v12); \draw (v12)--(v4);

\draw (v5)--(v13); \draw (v13)--(v15); \draw[opacity=0.4] (v15)--(v7);
\draw[opacity=0.4] (v6)--(v14); \draw[opacity=0.4] (v14)--(v16); \draw[opacity=0.4] (v16)--(v8);

\draw (v1)--(v17); \draw (v17)--(v19); \draw (v19)--(v2);
\draw (v3)--(v21); \draw[opacity=0.4] (v21)--(v23); \draw[opacity=0.4] (v23)--(v4);

\draw (v5)--(v18); \draw (v18)--(v20); \draw (v20)--(v6);
\draw[opacity=0.4] (v7)--(v22); \draw[opacity=0.4] (v22)--(v24); \draw[opacity=0.4] (v24)--(v8);

\draw (v9)--(v13); \draw (v13)--(v5);
\draw (v10)--(v14); \draw[opacity=0.4] (v14)--(v6);
\draw (v11)--(v15); \draw[opacity=0.4] (v15)--(v7);
\draw[opacity=0.4] (v12)--(v16); \draw[opacity=0.4] (v16)--(v8);

\draw (v9)--(v17); \draw (v13)--(v18);
\draw (v10)--(v19); \draw (v14)--(v20);
\draw (v11)--(v21); \draw[opacity=0.4] (v15)--(v22);
\draw[opacity=0.4] (v12)--(v23); \draw[opacity=0.4] (v16)--(v24);

\draw (v17)--(v18); \draw (v18)--(v20); \draw (v20)--(v19); \draw (v19)--(v17);
\draw[opacity=0.4] (v21)--(v22); \draw[opacity=0.4] (v22)--(v24); \draw[opacity=0.4] (v24)--(v23); \draw[opacity=0.4] (v23)--(v21);

\fill (v1) circle (0.5pt);
\fill (v2) circle (0.5pt);
\fill (v3) circle (0.5pt);
\fill (v4) circle (0.5pt);
\fill (v5) circle (0.5pt);
\fill (v6) circle (0.5pt);
\fill (v7) circle (0.5pt);
\fill (v8) circle (0.5pt);
\fill (v9) circle (0.5pt);
\fill (v10) circle (0.5pt);
\fill (v11) circle (0.5pt);
\fill (v12) circle (0.5pt);
\fill (v13) circle (0.5pt); 
\fill (v14) circle (0.5pt); 
\fill (v15) circle (0.5pt); 
\fill (v16) circle (0.5pt); 
\fill (v17) circle (0.5pt);
\fill (v18) circle (0.5pt); 
\fill (v19) circle (0.5pt);
\fill (v20) circle (0.5pt); 
\fill (v21) circle (0.5pt);
\fill (v22) circle (0.5pt); 
\fill (v23) circle (0.5pt);
\fill (v24) circle (0.5pt); 

\coordinate (pv13) at ( 0.414,  1, -1.514);
\coordinate (pv14) at ( 0.414, -1, -1.514);
\coordinate (pv15) at (-0.414,  1, -1.514);
\coordinate (pv16) at (-0.414, -1, -1.514);
\coordinate (pv18) at ( 1,  0.414, -1.514);
\coordinate (pv20) at ( 1, -0.414, -1.514);
\coordinate (pv22) at (-1,  0.414, -1.514);
\coordinate (pv24) at (-1, -0.414, -1.514);
\node (proj_label)[label={right:{\LARGE$\operatorname{proj}(\mathcal{X})$}}] at (+0.414,  1, -1.514){};
\node[above right = 200pt and 15pt of proj_label, label={right:{\LARGE$\mathcal{X}$}}]{};

\draw[fill=grey,opacity=0.2] (pv13)--(pv15)--(pv22)--(pv24)--(pv16)--(pv14)--(pv20)--(pv18)--cycle;

\draw[vector guide] (v13)--(pv13);
\draw[vector guide] (v14)--(pv14);
\draw[vector guide] (v15)--(pv15);
\draw[vector guide] (v16)--(pv16);
\draw[vector guide] (v18)--(pv18);
\draw[vector guide] (v20)--(pv20);
\draw[vector guide] (v22)--(pv22);
\draw[vector guide] (v24)--(pv24);
}

\end{tikzpicture}

    \vspace{-30pt}
    \caption{Example of the projection $\operatorname{proj}(\mathcal{X})$ of a three dimensional polytope $\mathcal{X}$ on the $xy$-plane.}
    \label{fig:poly_proj}
\end{figure}

The resulting constraint set is a convex and bounded polyhedron, i.e. a polytope. The set projection of $\mathcal{X}\subseteq\mathcal{Y}\times \mathcal{Z}$ onto $\mathcal{Z}$ is defined as in \cite{rakovic_reachability_2006}
\begin{equation}
    \operatorname{proj}_\mathcal{Z}(\mathcal{X})\doteq\{z\in\mathcal{Z}\ \vert\ \exists y\in\mathcal{Y}\text{ with } (z,y)\in\mathcal{X}\}\subseteq \mathcal{Z}.
\end{equation} 
This gives the values for the coupling variables  achievable without violating any of the constraints of the original set. Put differently, we can compute the FOR $\mathcal{F}_i$ as 
\begin{equation}\label{eq:prel_proj_ADP}
    \mathcal{F}_i \doteq \operatorname{proj}_{\mathcal{Z}_i}\mathcal{X}_i.
\end{equation}
An example of the projection of a three-dimensional polytope onto the xy-plane is shown in Figure~\ref{fig:poly_proj}.

\section{Complexity Reduction using Constrained Zonotopes}\label{sec:meth_cZ}
Different set representations offer different trade-offs in terms of representational accuracy, compactness, and suitability for certain operations, such as set projection, linear maps, and intersection. Therefore, the conversion into an appropriate set representation is critical to perform projection, i.e. aggregation, in high-dimensional problems such as multi-stage OPF.

\subsection{Approximation of Feasible Operating Region}
Given an affine set representation of the multi-stage OPF, the feasible set $\mathcal{X}_i$ is a convex and bounded polytope in half-space representation (H-rep) and can be expressed as
\begin{equation*} 
    \mathcal{P} \doteq \biggl\{ x\in \mathbb{R}^n \ \bigg| \ Ax\leq b \biggr\}, 
\end{equation*}
with a constraint matrix $A\in \mathbb{R}^{m\times n}$ and a constraint vector $b\in \mathbb{R}^m$.
This formulation enables the projection of the feasible set onto coupling variables using methods from computational geometry, such as Fourier–Motzkin elimination or vertex enumeration \cite[pp.~2-3]{jones_equality_nodate}. However, for multi-stage flexibility aggregation with storage units, the increased problem dimension renders these algorithms computationally intractable.
 
\subsection{Reformulation of the Feasible Set}
Zonotopes are well suited for projections, also in high dimensions. They are defined as
\begin{equation*}
    \mathcal{Z} \doteq \biggl\{ c + \sum_{i=1}^{n_g} \alpha_i G_{(\cdot,i)} \ \bigg| \ \alpha_i \in [-1,1]\biggr\},
\end{equation*}
with a center vector $c\in \mathbb{R}^n$, a generator matrix $G\in \mathbb{R}^{n\times n_g}$ and the scalar internal factors $\alpha_i$. Zonotope projection is performed via a simple sparse linear map and therefore can be computed very efficiently. In the following, the shorthand notation $\mathcal{Z}=\langle c,G\rangle_\mathcal{Z}$ is used for zonotopes. However, zonotopes are limited to model centrally symmetric convex polytopes.
To address this limitation Constrained Zonotopes ($\mathcal{CZ}\mathrm{s}$) extend zonotopes by adding linear constraints to the internal factors $\alpha_i$ \cite{CZ}. This results in the definition
\begin{equation*}\label{eq:prel_conZono}
\mathcal{CZ} \doteq \biggl\{c+\sum_{i=1}^{n_g}\alpha_i G_{(\cdot,i)} \ \bigg| \ \sum_{i=1}^{n_g}\alpha_i A_{(\cdot,i)}= b, \quad \alpha_i\in[-1,1]\biggl\},
\end{equation*}
with center vector $c\in \mathbb{R}^n$, generator matrix $G\in \mathbb{R}^{n\times n_g}$, constraint matrix $A\in \mathbb{R}^{m\times n_g}$, and constraint vector $b\in \mathbb{R}^{m}$. In the following, the shorthand notation $\mathcal{CZ}=\langle c,G,A,b\rangle_\mathcal{CZ}$ is used for constrained zonotopes.
Without limitations on the number of generators $n_g$ and the number of constraints $m$, $\mathcal{CZ}\mathrm{s}$ can express arbitrary convex polytopes and offer many of the computational advantages of zonotopes, like projection via matrix multiplication \cite{CZ}.

The conversion of the feasible set $\mathcal{X}_i$ of the multi-stage OPF into a $\mathcal{CZ}$ can be decomposed into three sequential stages, where the first two steps can be considered as an offline computation to achieve the desired set representation as $\mathcal{CZ}$ whereas the third step is set projection operation itself.

\subsubsection*{Determine a bounding zonotope}\label{sec:meth_bounding_zono}
Given a feasible set as convex and bounded polytope one can determine upper and lower bounds for each variable by solving two Linear Programms (LPs) for each variable $x_{i,k}$. This allows to determine the vectors of upper and lower variable bounds $ub_x$ and $lb_x$, which define the bounding zonotope $Z_b$ with center $c$ and generator matrix $G$ that over-approximates the given polytope $Z_b = \langle c,G\rangle_Z\subseteq \mathcal{P}$ as
\begin{subequations}
\label{eq:meth_bounding_zono}
\begin{align}
    c &= \tfrac{1}{2}(ub_x + lb_x)\\
    G &= \text{diag}(\tfrac{1}{2}(ub_x - lb_x))\\
    Z_b &= \langle c,\ G\rangle.
\end{align}
\end{subequations}
However, as shown in Subsection~\ref{subsec:eval_medium_grid} selecting the variable bounds sufficiently large is adequate without loss of accuracy, if for each variable $x_{i,k}$ the approximated variable bounds $\tilde{ub_x},\tilde{lb_x}$ exceed the variable bounds $\tilde{ub_x},\tilde{lb_x}$ with
\begin{equation}\label{eq:variable_bounds}
    \tilde{lb}_{x_i,k}\leq lb_{x_i,k}\leq x_{i,k}\leq ub_{x_i,k}\leq \tilde{ub}_{x_i,k}.
\end{equation}
While significantly reducing computation times in high dimensional problems such as multi-stage OPF, this may result in a highly over-approximative $Z_b$. However, the following step of the conversion into a $\mathcal{CZ}$ applies all constraints defining the feasible set $\mathcal{X}_i$ to the internal factors $\alpha$ of the $\mathcal{CZ}$ representation, i.e. effectively constraining the symmetric $Z_b$ to model the feasible set exactly. 


\subsubsection*{Adding constraints to the bounding zonotope}\label{subsec:meth_cZ_intersec}
Given the feasible set $\mathcal{X}_i$ of the multi-stage OPF as polytope $\mathcal{P}$ in H-rep, it is defined as an intersection of a finite number of halfspaces $H_i= \{ x \in \mathbb{R}^n \mid h^\top x \leq \zeta \}$. To reconstruct the polytope as $\mathcal{CZ}$, the bounding zonotope is successively intersected with each halfspace $H_i$. Based on the method in \cite{CZ_int}, each intersection updates the $\mathcal{CZ}$ as
\begin{equation}\label{eq:meth_cz_hs_int}
     \mathcal{CZ}_H = \langle c,\begin{bmatrix}
        G & 0
    \end{bmatrix}, \tilde{A}, \tilde{b}
    \rangle_{\mathcal{CZ}},
\end{equation}
with
\begin{equation}\label{eq:meth_cz_hs_int_update}
     \tilde{A}=\begin{bmatrix}
        A & 0 \\
        h^\top G & \frac{d_m}{2} 
    \end{bmatrix},\quad \tilde{b}=\begin{bmatrix}
        b \\
        \zeta-h^\top c - \frac{d_m}{2}
    \end{bmatrix},
\end{equation}
where $d_m$ is given by
\begin{equation}\label{eq:meth_cz_hs_int_dm}
    d_m = \zeta-h^\top c+ \sum_{i=1}^{n_g} | h^\top g_i|.
\end{equation}
This step assumes that the provided polytope is irredundant, so that no feasibility checks are required before intersection. Otherwise, it is beneficial to keep the representation compact by solving an LP to check whether a halfspace intersects the current $\mathcal{CZ}$ and therefore has to be added to the representation \cite{CZ_int}. The proposed approach for reformulation of the feasible set as $\mathcal{CZ}$ is summarized in the offline step in Algorithm~\ref{alg:meth_threaded_cz_int}.


\begin{algorithm}[t!]
\caption{Complexity reduction using $\mathcal{CZ}\mathrm{s}$.}
\label{alg:meth_threaded_cz_int}
\textbf{Require} Bounding zonotope \( \mathcal{Z} = \langle c, G \rangle \), halfspaces \( \mathcal{H} = \{ H_1, \dots, H_m \} \), defining polytope $P$\\
\textbf{Ensure} \( \mathcal{CZ} = \langle c, G, A, b \rangle \)\\
\textbf{Initialize} \( \mathcal{CZ} = \langle c, G, A, b \rangle \) with empty $A, b$\\
\text{Offline computation:} \begin{algorithmic}[1]
    \ForAll{\( H_i\) modeling static constraints \eqref{eq:PFeq}, \eqref{eq:netPwr}, and \eqref{eq:PF_df}-\eqref{eq:prel_storage_power_limits}}
        \State Compute \eqref{eq:meth_cz_hs_int_dm}
        \State Update \(\mathcal{CZ}\) as described in \eqref{eq:meth_cz_hs_int} and \eqref{eq:meth_cz_hs_int_update}
    \EndFor
\end{algorithmic}
\text{Online computation:} \begin{algorithmic}[1]
    \ForAll{\( H_i \) modeling dynamic constraints like \eqref{eq:genCstrBds}, \eqref{eq:ice_cream_ldf}, and \eqref{eq:prel_storage_dynamics}}
        \State Compute \eqref{eq:meth_cz_hs_int_dm}
        \State Update \(\mathcal{CZ}\) as described in \eqref{eq:meth_cz_hs_int} and \eqref{eq:meth_cz_hs_int_update}
    \EndFor
    \State Compute \eqref{eq:proj_cz}
\end{algorithmic}
\textbf{Return} $\operatorname{proj}(\mathcal{CZ})$
\end{algorithm}

The addition of a zero column to $G$ in \eqref{eq:meth_cz_hs_int} does not affect the result of $h^\top G$. This allows to parallelize the calculation of the rows for $A,b$ of the $\mathcal{CZ}$, which greatly benefits the required computation time.

\subsubsection*{Projection of the resulting constrained zonotope}\label{subsec:meth_cZ_proj_appr}
Finally, once the $\mathcal{CZ}$ is constructed, projection onto the desired subspace is performed via the sparse linear map $M_p$, yielding
\begin{equation}\label{eq:proj_cz}
    M_p \otimes \mathcal{CZ} = \langle M_pc,\ M_pG,\ A,\ b \rangle.
\end{equation}
This projection operation utilizes the computational benefits of the $\mathcal{CZ}$ representation, enabling fast and scalable dimensionality reduction via linear map. Thus, resulting in a computationally efficient online step of Algorithm~\ref{alg:meth_threaded_cz_int}. The additional steps of the online computation for dynamic sets are discussed in Subsection~\ref{subsec:eval_medium_grid}.
 
Unlike polytope projection that obscures details of sub-grids through variable elimination, the $\mathcal{CZ}$ projection maintains potentially confidential grid information inside $A\textrm{ and }b$ of the equality constraints. In their raw form the constraints might reflect local asset characteristics. However, the offline-online separation of the presented method allows for the seamless integration of privacy-preserving transformations in the offline step, i.e. with a non-singular private transformation matrix $T$ \cite{mangasarian_privacy-preserving_2012}
\begin{equation*}
    TA\alpha=Tb,
\end{equation*}
which masks the underlying structural data without affecting the communicated FOR.

\section{Intersection with Additional Linear Constraints}\label{sec:additional_intersec}
For more dynamic use cases, it is essential to limit the computation cost required to consider changes in the feasible set $\mathcal{X}_i$ due to control inputs or changing external conditions. Given the inherent separation of the set reformulation from the projection in the presented method shown in Algorithm~\ref{alg:meth_threaded_cz_int}, most of the computation can be performed offline. 
Furthermore, we can utilize that $\mathcal{CZ}\mathrm{s}$ allow to calculate halfspace and hyperplane intersections. This enables an efficient online update strategy for increasing problem sizes, as the underlying topology of the problem remains identical per time-step. We emphasize that, for dynamic use cases, the offline computation incorporates only the strictly static constraints modeling static constraints \eqref{eq:PFeq}, \eqref{eq:netPwr}, and \eqref{eq:PF_df}-\eqref{eq:prel_storage_power_limits} in $\mathcal{X}_i$. The offline step itself converts $\mathcal{X}_i$ into an equivalent $\mathcal{CZ}_{off}$, as elaborated in the previous section. This offline conversion accounts for the majority of the computational burden.
During online operation the offline $\mathcal{CZ}_{off}$ serves as the baseline. Time-varying bounds, such as the integration of revised forecast data in \eqref{eq:genCstrBds}, changing power demands or updated setpoint constraints, are integrated by adding the respective constraints via intersection to a copy of $\mathcal{CZ}_{off}$. This results in a new $\mathcal{CZ}_{on}$ modeling the feasible set for each iteration. Because $\mathcal{CZ}_{off}$ deliberately omits previous forecast data, there is no need to computationally remove or modify existing constraints. This avoids a repetition of the offline step entirely. The computational efficiency of adding constraints dynamically via intersections is evaluated in Subsection~\ref{sec:eval_additional_lin_cons}.

\section{Numerical Case Study} \label{sec:numerical_case_study} 
All numerical results in this section are performed on a system with an eight core and 16 thread Ryzen 9800X3D CPU paired with 64GB of 6000MHz RAM. Provided computation times are the average of ten executions. We use \texttt{JuMP.jl} \cite{lubin_jump_2023} and \texttt{MOSEK} \cite{mosek} to solve optimization problems. For convex sets, \texttt{Polyhedra.jl} allows to extract the set directly form the \texttt{JuMP.jl} model as polytope in H-representation.
For the problem formulation, we consider LinDistFlow with loss linearization as elaborated in Section~\ref{sec:probl_statement}. This allows to achieve a convex approximation while maintaining the crucial information of reactive power and voltage magnitude in the DSO grid model. As the underlying problem, we consider a multi-stage Optimal Power Flow (MS-OPF) problem, which allows flexibililty aggregation across multiple time steps and serves to evaluate computational scalability for horizons of up to $N=96$.

\subsection{Comparison of Computation Times}
First, we consider a small-scale, 4-bus radial distribution grid adapted from \texttt{MATPOWER}, called \texttt{case4dist} \cite{Matpower}. The grid has been extended to include two additional renewable generators at buses three and four, thereby increasing the system's overall active and reactive power flexibility. This small test case allows the comparison of the presented method with the baseline approach of polytope projection, which can still be evaluated over a short varying time horizon of the MS-OPF. While this projection-based approach does not scale well to larger grids, it remains suitable for evaluating performance in this compact setting.

\begin{table}[t!]
    \centering
    \caption{Computation times of the presented method and polytope projection on the 4-bus distribution grid model for up to $N=4$ timesteps. Fastest computation times in bold.}
    \label{tab:res_small_case}
    \setlength{\tabcolsep}{3pt}W
    \sisetup{
      scientific-notation=true,
      round-mode=places,
      round-precision=2,
      retain-zero-exponent=true,
      table-format=1.2e+1,  
      detect-weight=true,
      detect-inline-weight=math,
      exponent-product = \cdot
    }
    \begin{tabular}{ll
                    S 
                    S 
                    S }
        \toprule
        \textbf{Horizon \( N \)} & \textbf{Approach} & 
        \textbf{Offline [s]} & \textbf{Online [s]} & \textbf{Total [s]} \\
        \midrule
        \multirow{2}{*}{1} 
            & $\mathcal{CZ}$ & {\bfseries {\tablenum{0.00065725}}} & {\bfseries {\tablenum{1.42999e-6}}} & {\bfseries {\tablenum{0.00065865}}} \\
            & Polytope & {} & 0.0082195e+0 & 0.0082195e+0 \\
        \midrule
        \multirow{2}{*}{2} 
            & $\mathcal{CZ}$ & {\bfseries{\tablenum{0.0029587399}}} & {\bfseries{\tablenum{7.8000e-7}}} & {\bfseries{\tablenum{0.00295952}}} \\
            & Polytope & {} & 0.0581343e+0 & 0.0581343e+0 \\
        \midrule
        \multirow{2}{*}{3} 
            & $\mathcal{CZ}$ & {\bfseries{\tablenum{0.00728023}}} & {\bfseries{\tablenum{1.02e-6}}} & {\bfseries{\tablenum{0.00728125}}} \\
            & Polytope & {} & 0.2127586e+0 & 0.2127586e+0 \\
        \midrule
        \multirow{2}{*}{4} 
            & $\mathcal{CZ}$ & {\bfseries{\tablenum{0.011396410000000001}}} & {\bfseries{\tablenum{1.0299999999999999e-6}}} & {\bfseries{\tablenum{0.01139744}}} \\
            & Polytope & {} & 0.4627832e+0 & 0.4627832e+0 \\
        \bottomrule
    \end{tabular}
\end{table}

Table \ref{tab:res_small_case} shows the computation times for the multi-stage OPF on the 4-bus distribution grid model for up to $N=4$ timesteps. The results highlight that for small problem sizes, the presented method achieves significantly faster computation times, especially for the online computation of the set projection onto the coupling variables of active and reactive power at point of interconnection. While both methods model the convex feasible set of the multi-stage OPF based on the LinDistFlow formulation with loss linearization, the presented method achieves faster total computation times for all time horizon lengths depicted in Table \ref{tab:res_small_case} as the polytope projection, where no offline computation time is needed. This highlights, that for small problem sizes even the set conversion itself can be performed online, as its computation time is sufficiently fast.

\subsection{Evaluation on 15-bus Feeder}\label{subsec:eval_medium_grid}
Next, we analyze how the presented method performs on a medium-sized distribution grid model. 
\begin{figure}[t!]
    \vspace{-60pt}
    \centering
    \tikzset{sin v source/.style={
  circle,
  draw,
  append after command={
    \pgfextra{
    \draw
      ($(\tikzlastnode.center)!0.5!(\tikzlastnode.west)$)
       arc[start angle=180,end angle=0,radius=0.425ex] 
      (\tikzlastnode.center)
       arc[start angle=180,end angle=360,radius=0.425ex]
      ($(\tikzlastnode.center)!0.5!(\tikzlastnode.east)$) 
    ;
    }
  },
  scale=1.5,
 }
}

\tikzset{%
wind turbine/.pic={
  \tikzset{path/.style={fill, draw=black!40, line join=round}}
  \begin{scope}[scale=0.13]
  \path [path] 
    (-.25,0) arc (180:360:.25 and .0625) -- (.0625,3) -- (-.0625,3) -- cycle;
  \foreach \i in {90, 210, 330}{
    \ifcase#1
    \or
      \path [path, shift=(90:3), rotate=\i] 
        (.5,-.1875) arc (270:90:.5 and .1875) arc (90:-90:1.5 and .1875);
    \or
      \path [path, shift=(90:3), rotate=\i] 
        (0,0.125) -- (2,0.125) -- (2,0) -- (0.5,-0.375) -- cycle;
    \or
      \path [path, shift=(90:3), rotate=\i]
        (0,-0.125) arc (180:0:1 and 0.125) -- ++(0,0.125) arc (0:180:1 and 0.25) -- cycle;
    \fi
  }
  \path [path] (0,3) circle [radius=.25];
  \end{scope}
}}

\begin{tikzpicture}[
  node/.style={circle, draw, minimum size=5mm, inner sep=0pt},
  circ/.style={circle, draw, minimum size=5mm, inner sep=1pt, font=\footnotesize},
  ground/.style={thick},
  every node/.style={font=\small},
]
\scalebox{0.6}{
    \node[label=above:{1}] (Bus1) at (0,0) {};
    \node[right = 2.5 of Bus1, label=above:{2}] (Bus2) {};
    \node[right = 2.5 of Bus2, label=above:{3}] (Bus3) {};
    \node[right = 2.5 of Bus3, label=above:{4}] (Bus4) {};
    \node[right = 2.6 of Bus4, label=above:{5}] (Bus5) {};
    \node[below left = 4 and 0.5 of Bus2, label={[xshift=-3.0pt, yshift=-10pt]6}] (Bus6) {};
    \node[below = 3.6 of Bus1, label=above:{7}] (Bus7) {};
    \node[below = 2 of Bus6, label={[xshift=-3.0pt, yshift=-10pt]8}] (Bus8) {};
    \node[above left = 2 and 0.5 of Bus2, label={[xshift=-3.0pt, yshift=-10pt]9}] (Bus9) {};
    \node[above = 2 of Bus9, label={[xshift=-3.0pt, yshift=-10pt]10}] (Bus10) {};
    \node[right = 2.5 of Bus6, label={[xshift=-3.0pt, yshift=-10pt]11}] (Bus11) {};
    \node[below = 2 of Bus11, label={[xshift=-3.0pt, yshift=-10pt]12}] (Bus12) {};
    \node[below = 2 of Bus12, label={[xshift=-3.0pt, yshift=-10pt]13}] (Bus13) {};
    \node[above left = 2 and 0.5 of Bus4, label={[xshift=-3.0pt, yshift=-10pt]14}] (Bus14) {};
    \node[right = 2.5 of Bus11, label={[xshift=-3.0pt, yshift=-10pt]15}] (Bus15) {};
    \draw node[sacdcshape, anchor=south,scale=0.5] at (5.9, 0) (batt_conv) {};
    \draw node[twoportshape,above=0.3 of batt_conv.north,anchor=south,scale=0.5,fill={rgb,255:red,132; green,184; blue,25}] (bess_box) {};
    \draw node[battery1shape,above=0 of bess_box.center,anchor=center,scale=0.5,fill={rgb,255:red,132; green,184; blue,25}] (bess) {};
    
    \begin{scope}[every path/.style={draw=black, line width=2.0pt}]
        \draw (Bus1.center)-- ++(0,-2); 
        \draw (Bus2.center)-- ++(0,-2);
        \draw (Bus3.center)-- ++(0,-2); 
        \draw (Bus4.center)-- ++(0,-2); 
        \draw (Bus5.center)-- ++(0,-2);
        \draw (Bus6.east)-- ++(2,0);
        \draw (Bus7.center)-- ++(0,-2);
        \draw (Bus8.east)-- ++(2,0); 
        \draw (Bus9.east)-- ++(2,0); 
        \draw (Bus10.east)-- ++(2,0); 
        \draw (Bus11.east)-- ++(2,0); 
        \draw (Bus12.east)-- ++(2,0); 
        \draw (Bus13.east)-- ++(2,0); 
        \draw (Bus14.east)-- ++(2,0); 
        \draw (Bus15.east)-- ++(2,0); 
    \end{scope}
    
    \draw node[gridnode] (ext_grid) at (-1, -1) {};
    \draw (3.75,4) node [sin v source] (v1) {};
    \draw (3.75,-7.) node [sin v source] (v2) {};
    \draw (6.5,-9.225) node [sin v source] (v3) {};
    \path (4.3,3.75) pic {wind turbine=1};
    \path (4.3,-7.35) pic {wind turbine=1};
    \path (7.05,-9.575) pic {wind turbine=1};
    
    \begin{scope}[every path/.style={draw=black, line width=1.0pt}, every node/.style={draw,circle, minimum size = 23pt}]
        \draw (0,-1)-- ++(11.1,0)-- ++(-9.7,0) node[yshift=14pt, draw, circle]{1}-- ++(2.75,0) node[yshift=14pt, draw, circle]{2}-- ++(2.75,0) node[yshift=14pt, draw, circle]{3}-- ++(2.875,0) node[yshift=14pt, draw, circle]{4}; 
        \draw (2.75,-1.5)-- ++(0.4,0)-- ++(0,-5)-- ++(0,3.5) node[xshift=-14pt, draw, circle]{7}-- ++(0,-2.4) node[xshift=-14pt, draw, circle]{9};
        \draw (5.5,-1.5)-- ++(0.4,0)-- ++(0,-7.25)-- ++(0,1.125) node[xshift=-14pt, draw, circle]{12}-- ++(0,2.25) node[xshift=-14pt, draw, circle]{11}-- ++(0,2.25) node[xshift=-14pt, draw, circle]{10};
        \draw (8.25,-1.5)-- ++(0.4,0)-- ++(0,-2.75)-- ++(0,1.125) node[xshift=-14pt, draw, circle]{14};
        
        \draw (2.75,-0.5)-- ++(0.4,0)-- ++(0,1.75) node[xshift=-14pt, draw, circle]{5}-- ++(0,2.125) node[xshift=-14pt, draw, circle]{6};
        \draw (2.75,-0.5)-- ++(0.4,0)-- ++(0,5);
        \draw (8.3,-0.5)-- ++(0.35,0)-- ++(0,2.75)-- ++(0,-1) node[xshift=-14pt, draw, circle]{13};
    
        \draw (2.5, -4.25)-- ++(0,-0.625)-- ++(-2.5,0)-- ++(1.375,0) node[yshift=-14pt, draw, circle]{8};
        \draw[-stealth](3.75,2.25) --++ (0,-0.5);
        \draw[-stealth](3.75,-4.25) --++ (0,-0.5);
    
        \draw[-stealth](6.5,-4.25) --++ (0,-0.5);
        \draw[-stealth](6.5,-6.5) --++ (0,-0.5);
    
        \draw[-stealth](9.25,2.25) --++ (0,-0.5);
        \draw[-stealth](9.25,-4.25) --++ (0,-0.5);
    
        \draw[-stealth] (2.75,-1.5)-- ++(-0.4,0)--++ (0,-0.5);
        \draw[-stealth] (5.5,-1.5)-- ++(-0.4,0)--++ (0,-0.5);
        \draw[-stealth] (8.25,-1.5)-- ++(-0.4,0)--++ (0,-0.5);
        \draw[-stealth] (11.1,-1.5)-- ++(0.4,0)--++ (0,-0.5);
        \draw[-stealth] (0, -5.375)-- ++(-0.4,0)--++ (0,-0.5);
        
        \draw (ext_grid.east)-- ++(1,0);
        \draw (v1)-- ++(0,0.5);
        \draw (v2)-- ++(0,0.5);
        \draw (v3)-- ++(0,0.5);
    
         \draw (5.5,-0.5)-- ++(0.4,0)--++(0,0.5);
         \draw (batt_conv.north)-- (bess_box.south);
    
    \end{scope}
}

\end{tikzpicture}
    \vspace{-100pt}
    \caption{Topology of the 15-bus radial distribution grid with battery highlighted in green.}
    \label{fig:grid_model}
\end{figure}
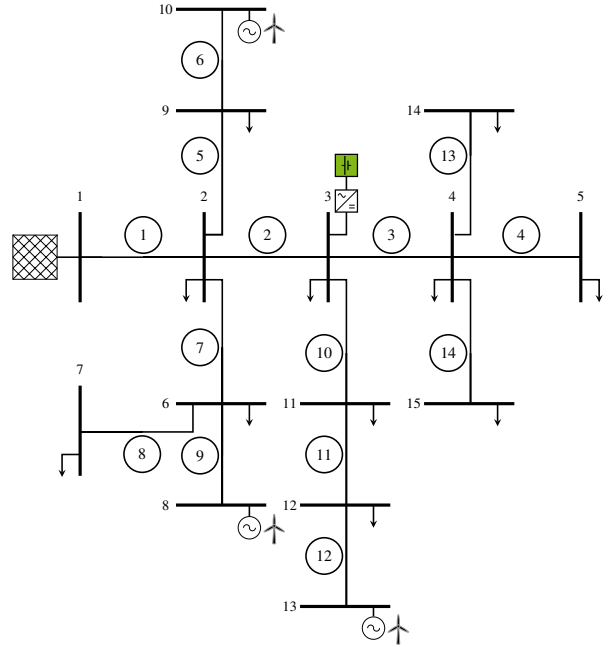
\begin{figure}[t!]
    \centering
    \includesvg[width=1\linewidth]{problem_size_scaling.svg}
    \caption{Scaling of dimensionality and number of constraints of the 15-bus grid model with time horizon length $N$.}
    \label{fig:eval_problem_size}
\end{figure}
Therefore, we consider the 15-bus grid model \texttt{case15nbr} from \cite{Matpower,case_15}, where the loads on buses
8,10 and 13 are replaced by renewable generators with double the active power limit and an added storage unit at bus 3 with a capacity of 1~\si{MWh}, active power of 1~\si{MW}, and an initial SOC of 100~\si{\percent}. The adapted grid is depicted in Figure~\ref{fig:grid_model}.  Figure~\ref{fig:eval_problem_size} illustrates the scaling of the optimization problem as the time horizon $N$ increases. While the underlying 15-bus grid model can be considered as a rather small MV network, the scaling of the problem size with increasing time horizon length highlights the applicability of the approach to large-scale problems.

\subsubsection*{Approximation quality of reformulation based approaches}\label{sec:res_area_cz_lp_scaling}
\begin{figure}[t!]
    \centering
    \includesvg[width=.9\linewidth]{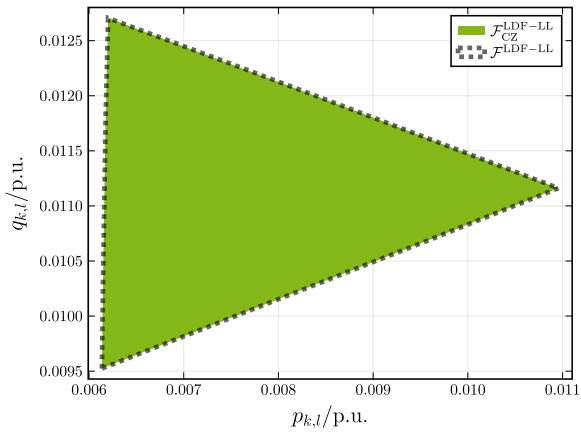}
    \caption{FOR at the slack bus modeled as polytope and as $\mathcal{CZ}$ using sufficiently large variable bounds.}
    \label{fig:lp_calc_cz_comp_area}
\end{figure}
\begin{table}[t!]
    \centering
    \caption{Computation times of the presented method on the 15-bus distribution grid model for up to $N=96$ timesteps.}
    \label{tab:res_case_15}
    
    \sisetup{
      scientific-notation=true,
      round-mode=places,
      round-precision=2,
      retain-zero-exponent=true,
      table-format=1.2e+1,  
      detect-weight=true,
      detect-inline-weight=math,
      exponent-product = \cdot
    }
    \begin{tabular}{l
                    S 
                    S 
                    S }
        \toprule
        \textbf{Horizon \( N \)} & \textbf{Offline [s]} & \textbf{Online [s]} & \textbf{Total [s]} \\
        \midrule
        \multirow{1}{*}{12} 
            & 1.17887226 & 2.9520000000000002e-5 & 1.17890178 \\
        \midrule
        \multirow{1}{*}{24} 
            & 5.48808776 & 0.00011395 & 5.48820171 \\
        \midrule
        \multirow{1}{*}{36} 
            & 23.795031979999997 & 0.00011068999999999999 & 23.79514267 \\
        \midrule
        \multirow{1}{*}{48} 
            & 61.923130689999994 & 0.00014536999999999997 & 61.92327606 \\
        \midrule
        \multirow{1}{*}{60} 
            & 124.31393603000001 & 0.00150408 & 124.31544011 \\
        \midrule
        \multirow{1}{*}{72} 
            & 267.92545057 & 0.0032322 & 267.92868277 \\
        \midrule
        \multirow{1}{*}{84} 
            & 311.89421432000006 & 0.0022251899999999993 & 311.89643951 \\
        \midrule
        \multirow{1}{*}{96} 
            & 535.23161739 & 0.00284774 & 535.23446513 \\
        \bottomrule
    \end{tabular}
\end{table}
 First, we demonstrate the equivalence of the reformulation of the feasible set as $\mathcal{CZ}$. Figure~\ref{fig:lp_calc_cz_comp_area} shows that approximating the variable bounds does not impact the accuracy of the resulting $\mathcal{CZ}$. The FOR modeled as $\mathcal{CZ}$ $\mathcal{F}_{\mathcal{CZ}}^{LDF-LL}$ is equivalent to the respective FOR modeled as polytope $\mathcal{F}^{LDF-LL}$. As described in Subsection \ref{sec:meth_bounding_zono}, this holds as long as the bounding zonotope over-approximates the polytope and all constraints defining the polytope are added to the $\mathcal{CZ}$ representation in the next step of the conversion. Hence, all following computations use sufficiently large variable bounds to decrease computation time. For a comparison of the approximation quality of $\mathcal{F}^{LDF-LL}$ with the true flexibility region we refer to \cite{ADP}.

Table \ref{tab:res_case_15} shows the respective timeseries of required computation times for the reformulation of the feasible set as $\mathcal{CZ}$. For the given results, the constraint intersection is performed in parallel. The results clearly show that despite the fact that the computation time scales noticeably, it remains suitable for offline computation on the used test problem, even for a considered time horizon length of $N=96$. Further, $\mathcal{CZ}$ projection can be performed in sub-second times on the 15-bus distribution grid case, even for large time horizons $N$. This highlight the suitability for the presented method to reduce complexity in grid operation via aggregation. Paired with acceptable conversion computation times, this allows the application of aggregation to considerably larger problem sizes, by using the presented method to model the feasible sets of the sub problems as $\mathcal{CZ}$.

\subsubsection*{Computation times for set modification}\label{sec:eval_additional_lin_cons} 
To reduce the repeated necessity of performing the set conversion due to small changes in the feasible set, e.g. set point alterations, this section briefly evaluates the modification of an existing $\mathcal{CZ}$ with additional halfspace or hyperplane intersection as discussed in Section~\ref{sec:additional_intersec}. 


\begin{table}[t!]
    \centering
    \caption{Computation times required to add an additional constraint to an existing $\mathcal{CZ}$ for up to $N=96$ timesteps.}
    \label{tab:res_big_case}
    
    \sisetup{
      scientific-notation=true,
      round-mode=places,
      round-precision=2,
      retain-zero-exponent=true,
      table-format=1.2e+1,  
      detect-weight=true,
      detect-inline-weight=math,
      exponent-product = \cdot
    }
    \begin{tabular}{l
                    S 
                    }
        \toprule
        \textbf{Horizon \( N \)} & \textbf{Computation time [s]} \\
        \midrule
        \multirow{1}{*}{12} 
            & 0.03347493 \\
        \midrule
        \multirow{1}{*}{24} 
            & 0.24586901999999994 \\
        \midrule
        \multirow{1}{*}{36} 
            & 0.25688615000000004\\
        \midrule
        \multirow{1}{*}{48} 
            & 0.44244116\\
        \midrule
        \multirow{1}{*}{60} 
            & 1.3807793\\
        \midrule
        \multirow{1}{*}{72} 
            & 3.69932956\\
        \midrule
        \multirow{1}{*}{84} 
            & 4.82410281\\
        \midrule
        \multirow{1}{*}{96} 
            & 6.5599988499999995\\
        \bottomrule
    \end{tabular}
\end{table}

Table~\ref{tab:res_big_case} depicts the time necessary to add an individual affine constraint to a $\mathcal{CZ}$ modeling the multi-stage OPF for increasing time horizon lengths $N$. The results show, that adding additional affine constraints is in fact viable as the required computation time remains sufficiently low, i.e. far below the considered time step size of fifteen minutes, even for large problem size. 

These result highlight the potential use of these set representations for dynamic sets. For instance, this can be used in ADP to model a feasible set for a fixed set of constraints and add more dynamic constraints, e.g. time-varying generator setpoints, via further intersection as necessary.

\subsubsection*{Time-dependent feasible operating region}
\begin{figure}[t!]
    \centering
    \includesvg[width=.8\linewidth]{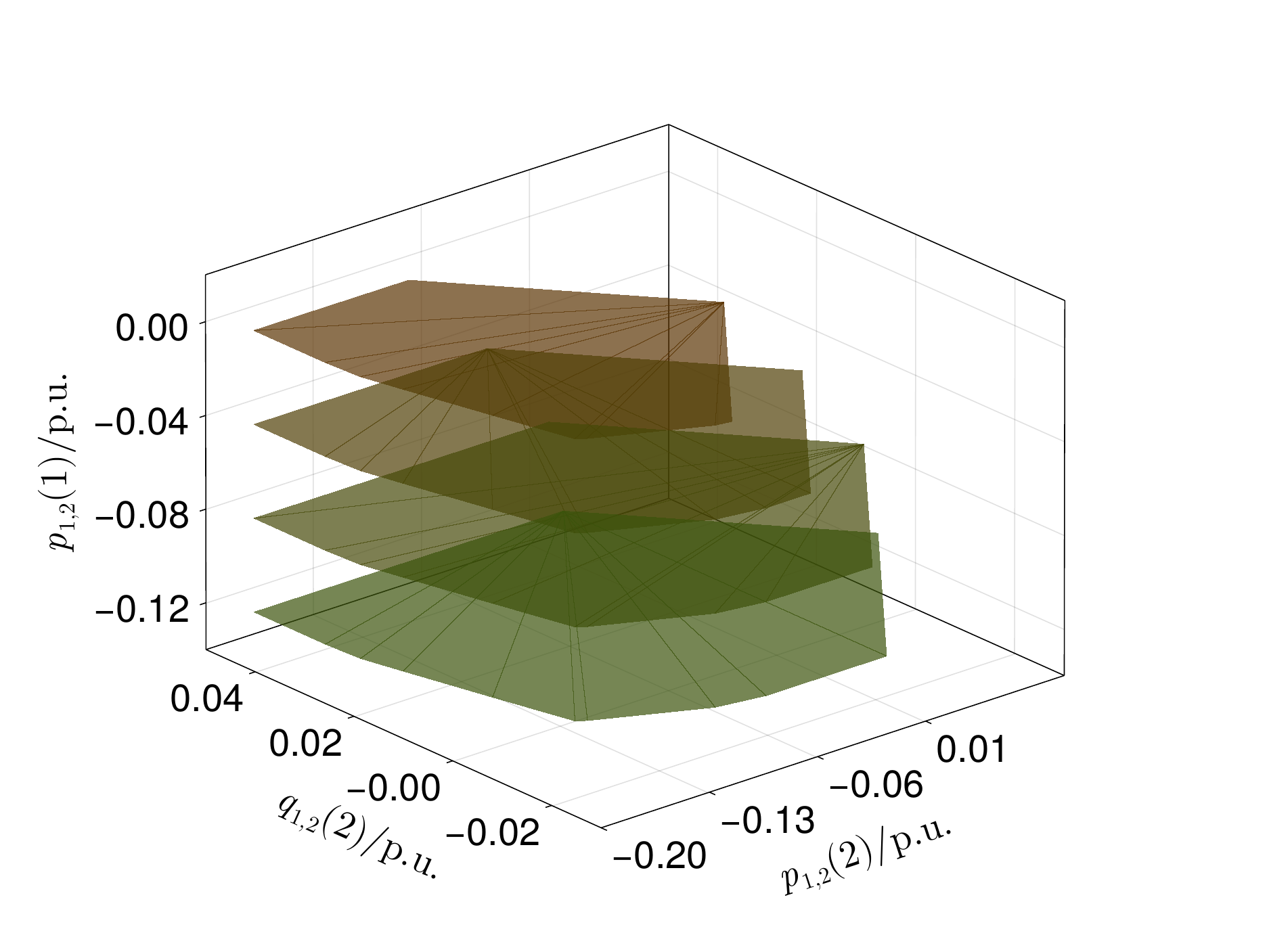}
    \caption{Dependency of $p_{1,2}(1)$ FOR of active and reactive power at the following timestep.}
    \label{fig:time_dependent_for}
\end{figure}

Given the significantly reduced computation times, we can now consider a MS-OPF problem on the 15-bus distribution grid case and therefore the FOR at point of interconnection for the coupling variables of active and reactive power along the time horizon. Figure~\ref{fig:time_dependent_for} shows the available active and reactive power values of the multidimensional FOR at time step $k=2$, denoted by $p_{1,2}(2)$ and $q_{1,2}(2)$, depending on the provided active power at the point of interconnection at $k=1$, $p_{1,2}(1)$. One can see that the $p_{1,2}(2)$/$q_{1,2}(2)$ area becomes significantly smaller with increasing $p_{1,2}(1)$. This is due to the battery dynamics \eqref{eq:prel_storage_dynamics}, i.e. the available energy in the battery is lowered for timestep $k=2$ depending on the power drawn in timestep $k=1$. Otherwise, the available flexibility remains consistent as the loads as well as generation units in the 15-bus distribution grid model are considered to be static over time.

This highlights the necessity to consider the multi-stage OPF for the DSO subproblems in order to ensure a feasible solution for consecutive timesteps and to improve flexibilty utilization.

\section{Conclusion and Future Work} \label{sec:conclusion}
This paper addresses the central computational bottleneck arising in projection-based determination of the FOR for affine grid models via tools of computational geometry. By utilizing $\mathcal{CZ}\mathrm{s}$ as beneficial set representations we can significantly reduce the computation times required for flexibility aggregation via set projection of convex feasible sets while achieving the same accuracy of original polytope based approaches. Hence, this approach combines the computational efficiency of zonotope based approximations for set projection with the capability to model non-symmetric polytopes. The computational efficiency of the presented method allows to perform time-dependent flexibility aggregation in the setting of multi-stage OPF and on large-scale problems without loss of accuracy. Further we show that most of the computation can be performed offline and set modifications can be considered efficiently.

Future work will focus on further application in the proposed use cases and a more detailed numerical comparison of the the representational accuracy to alternative approaches like original zonotope projection. Additionally future work will focus on the extension of the presented method to non-convex OPF formulations as well as the charging-discharging efficiency of batteries and their nonlinear dynamics.

\AtNextBibliography{\small}
\printbibliography

@article{Engelmann2025,
  title = {Approximate {{Dynamic Programming With Feasibility Guarantees}}},
  author = {Engelmann, Alexander and Bandeira, Ma{\'i}sa Beraldo and Faulwasser, Timm},
  year = {2025},
  journal = {IEEE Transactions on Control of Network Systems},
  pages = {1565 -- 1576},
  issn = {2325-5870},
  doi = {10.1109/TCNS.2025.3526715},
  urldate = {2025-01-15}
}

@article{capitanescu2018,
  title = {{{TSO}}\textendash{{DSO}} Interaction: {{Active}} Distribution Network Power Chart for {{TSO}} Ancillary Services Provision},
  shorttitle = {{{TSO}}\textendash{{DSO}} Interaction},
  author = {Capitanescu, Florin},
  year = {2018},
  month = oct,
  journal = {Electric Power Systems Research},
  volume = {163},
  pages = {226--230},
  issn = {0378-7796},
  doi = {10.1016/j.epsr.2018.06.009},
  langid = {english}
}

@article{nazir2022,
  title = {Grid-{{Aware Aggregation}} and {{Realtime Disaggregation}} of {{Distributed Energy Resources}} in {{Radial Networks}}},
  author = {Nazir, Nawaf and Almassalkhi, Mads},
  year = {2022},
  month = may,
  journal = {IEEE Transactions on Power Systems},
  volume = {37},
  number = {3},
  pages = {1706--1717},
  issn = {1558-0679},
  doi = {10.1109/TPWRS.2021.3121215}
}

@article{sarstedt2022,
  title = {Monetarization of the {{Feasible Operation Region}} of {{Active Distribution Grids Based}} on a {{Cost-Optimal Flexibility Disaggregation}}},
  author = {Sarstedt, Marcel and Hofmann, Lutz},
  year = {2022},
  journal = {IEEE Access},
  volume = {10},
  pages = {5402--5415},
  issn = {2169-3536},
  doi = {10.1109/ACCESS.2022.3140871}
}

@article{silva2018,
  title = {Estimating the Active and Reactive Power Flexibility Area at the {{TSO-DSO}} Interface},
  author = {Silva, Jo{\~a}o and Sumaili, Jean and Bessa, Ricardo J. and Seca, Lu{\'i}s and Matos, Manuel A. and Miranda, Vladimiro and Caujolle, Mathieu and Goncer, Bel{\'e}n and {Sebastian-Viana}, Maria},
  year = {2018},
  month = sep,
  journal = {IEEE Transactions on Power Systems},
  volume = {33},
  number = {5},
  pages = {4741--4750},
  issn = {1558-0679},
  doi = {10.1109/TPWRS.2018.2805765}
}

@article{baran1989-0,
  title = {Optimal Sizing of Capacitors Placed on a Radial Distribution System},
  author = {Baran, M. and Wu, F.F.},
  year = {1989},
  month = jan,
  journal = {IEEE Transactions on Power Delivery},
  volume = {4},
  number = {1},
  pages = {735--743},
  issn = {1937-4208},
  doi = {10.1109/61.19266}
}

@inproceedings{contreras2019,
  title = {Time-{{Based Aggregation}} of {{Flexibility}} at the {{TSO-DSO Interconnection Point}}},
  booktitle = {{{IEEE Power}} \& {{Energy Society General Meeting}} ({{PESGM}})},
  author = {Contreras, Daniel A. and Rudion, Krzysztof},
  year = {2019},
  month = aug,
  pages = {1--5},
  issn = {1944-9933},
  doi = {10.1109/PESGM40551.2019.8973421}
}

@article{sarstedt2021,
  title = {Survey and {{Comparison}} of {{Optimization-Based Aggregation Methods}} for the {{Determination}} of the {{Flexibility Potentials}} at {{Vertical System Interconnections}}},
  author = {Sarstedt, Marcel and Klu{\ss}, Leonard and Gerster, Johannes and Meldau, Tobias and Hofmann, Lutz},
  year = {2021},
  month = jan,
  journal = {Energies},
  volume = {14},
  number = {3},
  pages = {687},
  publisher = {{Multidisciplinary Digital Publishing Institute}},
  issn = {1996-1073},
  doi = {10.3390/en14030687},
  urldate = {2023-07-26},
  copyright = {http://creativecommons.org/licenses/by/3.0/},
  langid = {english}
}

@misc{Bandeira2025,
  title = {Complexity {{Reduction}} for {{TSO-DSO Coordination}}: {{Flexibility Aggregation}} vs. {{Distributed Optimization}}},
  shorttitle = {Complexity {{Reduction}} for {{TSO-DSO Coordination}}},
  author = {Bandeira, Ma{\'i}sa Beraldo and Engelmann, Alexander and Faulwasser, Timm},
  year = {2025},
  month = sep,
  number = {arXiv:2509.10595},
  eprint = {2509.10595},
  primaryclass = {eess},
  publisher = {arXiv},
  doi = {10.48550/arXiv.2509.10595},
  urldate = {2025-09-16},
  archiveprefix = {arXiv},
}

@article{ozturk2022,
  title = {Aggregation of {{Demand-Side Flexibilities}}: {{A Comparative Study}} of {{Approximation Algorithms}}},
  shorttitle = {Aggregation of {{Demand-Side Flexibilities}}},
  author = {{\"O}zt{\"u}rk, Emrah and Rheinberger, Klaus and Faulwasser, Timm and Worthmann, Karl and Prei{\ss}inger, Markus},
  year = {2022},
  month = jan,
  journal = {Energies},
  volume = {15},
  number = {7},
  pages = {2501},
  publisher = {Multidisciplinary Digital Publishing Institute},
  issn = {1996-1073},
  doi = {10.3390/en15072501},
  urldate = {2025-10-06},
  copyright = {http://creativecommons.org/licenses/by/3.0/},
  langid = {english}
}

@inproceedings{jha2020,
  title = {Coordinated {{Voltage Control}} for {{Conservation Voltage Reduction}} in {{Power Distribution Systems}}},
  booktitle = { {{IEEE Power}} \& {{Energy Society General Meeting}} ({{PESGM}})},
  author = {Jha, Rahul Ranjan and Dubey, Anamika},
  year = {2020},
  month = aug,
  pages = {1--5},
  issn = {1944-9933},
  doi = {10.1109/PESGM41954.2020.9282024}
}

@article{baran1989,
  title = {Network Reconfiguration in Distribution Systems for Loss Reduction and Load Balancing},
  author = {Baran, M.E. and Wu, F.F.},
  year = {1989},
  month = apr,
  journal = {IEEE Transactions on Power Delivery},
  volume = {4},
  number = {2},
  pages = {1401--1407},
  issn = {1937-4208},
  doi = {10.1109/61.25627}
}

@inproceedings{wang2023,
  title = {A {{Projection-Based Approach}} for {{Distributed Energy Resources Aggregation}}},
  booktitle = {2023 {{IEEE PES Innovative Smart Grid Technologies Europe}} ({{ISGT EUROPE}})},
  author = {Wang, Yiran and Zhong, Haiwang and Ruan, Guangchun},
  year = {2023},
  month = oct,
  pages = {1--5},
  doi = {10.1109/ISGTEUROPE56780.2023.10408587},
  urldate = {2025-10-13}
}

@article{ADP,
	title = {An {ADP} framework for flexibility and cost aggregation: Guarantees and open problems},
	volume = {234},
	issn = {0378-7796},
	url = {https://www.sciencedirect.com/science/article/pii/S0378779624007041},
	doi = {10.1016/j.epsr.2024.110818},
	shorttitle = {An {ADP} framework for flexibility and cost aggregation},
	note = {Art. no. 110818},
	journaltitle = {Electric Power Systems Research},
	shortjournal = {Electr. Power Syst. Res.},
	author = {Beraldo Bandeira, Maísa and Faulwasser, Timm and Engelmann, Alexander},
	urldate = {2025-07-27},
	date = {2024-09-01},
}

@article{case_15,
	title = {{DG} planning with amalgamation of economic and reliability considerations},
	volume = {73},
	issn = {0142-0615},
	url = {https://www.sciencedirect.com/science/article/pii/S0142061515002112},
	doi = {10.1016/j.ijepes.2015.05.006},
	pages = {273--282},
	journaltitle = {International Journal of Electrical Power \& Energy Systems},
	shortjournal = {Int. J. Elect. Power Eng. Syst.},
	author = {Battu, Neelakanteshwar Rao and Abhyankar, A. R. and Senroy, Nilanjan},
	urldate = {2025-07-27},
	date = {2015-12-01},
}

@article{MP_CC_Aggregation,
	title = {Barriers and Insights to Compute Multi-Period Cost Curves of Active Power Aggregated Flexibility from Distribution Systems for {TSO}-{DSO} Coordination},
	issn = {1558-0679},
	url = {https://ieeexplore.ieee.org/document/10982438},
	doi = {10.1109/TPWRS.2025.3566704},
	pages = {1--12},
	journaltitle = {{IEEE} Transactions on Power Systems},
	author = {Capitanescu, Florin},
	urldate = {2025-07-27},
	date = {2025},
}

@article{CZ_int,
	title = {Set operations and order reductions for constrained zonotopes},
	volume = {139},
	issn = {0005-1098},
	url = {https://www.sciencedirect.com/science/article/pii/S0005109822000498},
	doi = {10.1016/j.automatica.2022.110204},
	note = {Art. no. 110204},
	journaltitle = {Automatica},
	shortjournal = {Automatica},
	author = {Raghuraman, Vignesh and Koeln, Justin P.},
	urldate = {2025-07-27},
	date = {2022-05-01},
}

@article{CZ,
	title = {Constrained zonotopes: A new tool for set-based estimation and fault detection},
	volume = {69},
	issn = {0005-1098},
	url = {https://www.sciencedirect.com/science/article/pii/S0005109816300772},
	doi = {10.1016/j.automatica.2016.02.036},
	shorttitle = {Constrained zonotopes},
	pages = {126--136},
	journaltitle = {Automatica},
	shortjournal = {Automatica},
	author = {Scott, Joseph K. and Raimondo, Davide M. and Marseglia, Giuseppe Roberto and Braatz, Richard D.},
	urldate = {2025-07-27},
	date = {2016-07-01},
}

@article{Matpower,
	title = {{MATPOWER}: Steady-State Operations, Planning, and Analysis Tools for Power Systems Research and Education},
	volume = {26},
	issn = {1558-0679},
	url = {https://ieeexplore.ieee.org/document/5491276},
	doi = {10.1109/TPWRS.2010.2051168},
	shorttitle = {{MATPOWER}},
	pages = {12--19},
	number = {1},
	journaltitle = {{IEEE} Transactions on Power Systems},
	author = {Zimmerman, Ray Daniel and Murillo-Sánchez, Carlos Edmundo and Thomas, Robert John},
	urldate = {2025-07-27},
	date = {2011-02},
}

@article{rakovic_reachability_2006,
	title = {Reachability analysis of discrete-time systems with disturbances},
	volume = {51},
	issn = {1558-2523},
	url = {https://ieeexplore.ieee.org/document/1618830},
	doi = {10.1109/TAC.2006.872835},
	pages = {546--561},
	number = {4},
	journaltitle = {{IEEE} Transactions on Automatic Control},
	author = {Rakovic, S.V. and Kerrigan, E.C. and Mayne, D.Q. and Lygeros, J.},
	urldate = {2025-07-27},
	date = {2006-04},
}

@techreport{jones_equality_nodate,
	title = {Equality Set Projection: A new algorithm for the projection of polytopes in halfspace representation},
	url = {https://infoscience.epfl.ch/handle/20.500.14299/71951},
	author = {{C. Jones, E. Kerrigan, and J. Maciejowski}},
	date = {2004},
    institution = {Eng. Dept., Cambridge Univ.},
    location={Cambridge, {U.K.}},
    number = {TR.463}
}

@article{lubin_jump_2023,
	title = {{JuMP} 1.0: recent improvements to a modeling language for mathematical optimization},
	volume = {15},
	issn = {1867-2957},
	shorttitle = {{JuMP} 1.0},
	url = {https://doi.org/10.1007/s12532-023-00239-3},
	doi = {10.1007/s12532-023-00239-3},
	number = {3},
	urldate = {2025-07-31},
	journal = {Mathematical Programming Computation},
	author = {Lubin, Miles and Dowson, Oscar and Garcia, Joaquim Dias and Huchette, Joey and Legat, Benoît and Vielma, Juan Pablo},
	month = sep,
	year = {2023},
	pages = {581--589},
}

@article{mosek,
   author = "MOSEK",
   title = "The {MOSEK} Optimizer {API} for Julia manual. {Version} 11.0.",
    howpublished = {[Online]},
   year = 2025,
   note = {Accessed: Jul.~31,~2025. [Online]. Available: \url{https://docs.mosek.com/latest/juliaapi/index.html}}
 }

@ARTICLE{agg_disagg_zonotope,
  author={Müller, Fabian L. and Szabó, Jácint and Sundström, Olle and Lygeros, John},
  journal={IEEE Transactions on Smart Grid}, 
  title={Aggregation and Disaggregation of Energetic Flexibility From Distributed Energy Resources}, 
  year={2019},
  volume={10},
  number={2},
  pages={1205-1214},
  doi={10.1109/TSG.2017.2761439}}

@INPROCEEDINGS{zono_union,
  author={Nazir, Md Salman and Hiskens, Ian A. and Bernstein, Andrey and Dall'Anese, Emiliano},
  booktitle={2018 IEEE Conference on Decision and Control (CDC)}, 
  title={Inner Approximation of Minkowski Sums: A Union-Based Approach and Applications to Aggregated Energy Resources}, 
  year={2018},
  volume={},
  number={},
  pages={5708-5715},
  doi={10.1109/CDC.2018.8618731}}

@ARTICLE{Brett2016,
  author={Robbins, Brett A. and Domínguez-García, Alejandro D.},
  journal={IEEE Transactions on Power Systems}, 
  title={Optimal Reactive Power Dispatch for Voltage Regulation in Unbalanced Distribution Systems}, 
  year={2016},
  volume={31},
  number={4},
  pages={2903-2913},
  doi={10.1109/TPWRS.2015.2451519}}

@article{mangasarian_privacy-preserving_2012,
	title = {Privacy-preserving horizontally partitioned linear programs},
	volume = {6},
	copyright = {http://www.springer.com/tdm},
	issn = {1862-4472, 1862-4480},
	url = {http://link.springer.com/10.1007/s11590-010-0268-9},
	doi = {10.1007/s11590-010-0268-9},
	language = {en},
	number = {3},
	urldate = {2026-03-02},
	journal = {Optimization Letters},
	author = {Mangasarian, Olvi L.},
	month = mar,
	year = {2012},
	pages = {431--436},
}

@ARTICLE{Chen2020,
  author={Chen, Xin and Dall’Anese, Emiliano and Zhao, Changhong and Li, Na},
  journal={IEEE Transactions on Smart Grid}, 
  title={Aggregate Power Flexibility in Unbalanced Distribution Systems}, 
  year={2020},
  volume={11},
  number={1},
  pages={258-269},
  doi={10.1109/TSG.2019.2920991}}

@ARTICLE{Chen2021,
  author={Chen, Xin and Li, Na},
  journal={IEEE Transactions on Smart Grid}, 
  title={Leveraging Two-Stage Adaptive Robust Optimization for Power Flexibility Aggregation}, 
  year={2021},
  volume={12},
  number={5},
  pages={3954-3965},
  doi={10.1109/TSG.2021.3068341}}

@ARTICLE{Cui2021,
  author={Cui, Bai and Zamzam, Ahmed and Bernstein, Andrey},
  journal={IEEE Control Systems Letters}, 
  title={Network-Cognizant Time-Coupled Aggregate Flexibility of Distribution Systems Under Uncertainties}, 
  year={2021},
  volume={5},
  number={5},
  pages={1723-1728},
  doi={10.1109/LCSYS.2020.3045080}}

@ARTICLE{2Tan2024,
  author={Tan, Zhenfei and Yu, Ao and Zhong, Haiwang and Zhang, Xianfeng and Xia, Qing and Kang, Chongqing},
  journal={CSEE Journal of Power and Energy Systems}, 
  title={Optimal Virtual Battery Model for Aggregating Storage-Like Resources with Network Constraints}, 
  year={2024},
  volume={10},
  number={4},
  pages={1843-1847},
  doi={10.17775/CSEEJPES.2022.04090}}

@INPROCEEDINGS{Gonzales2018,
  author={Mayorga Gonzalez, D. and Hachenberger, J. and Hinker, J. and Rewald, F. and Häger, U. and Rehtanz, C. and Myrzik, J.},
  booktitle={2018 Power Systems Computation Conference (PSCC)}, 
  title={Determination of the Time-Dependent Flexibility of Active Distribution Networks to Control Their TSO-DSO Interconnection Power Flow}, 
  year={2018},
  volume={},
  number={},
  pages={1-8},
  doi={10.23919/PSCC.2018.8442865}}

@article{ozturk2024alleviating,
title={Alleviating the curse of dimensionality in minkowski sum approximations of storage flexibility},
author={{\"O}zt{\"u}rk, Emrah and Faulwasser, Timm and Worthmann, Karl and Prei{\ss}inger, Markus and Rheinberger, Klaus},
journal={IEEE Transactions on Smart Grid},
volume={15},
number={6},
pages={5733--5743},
year={2024},
publisher={IEEE}
}



\end{document}